\documentclass[11pt,reqno]{amsart}
\usepackage{amssymb,verbatim,enumerate,ifthen}
\usepackage[mathscr]{eucal}
\def\N{\mathbb{N}}
\def\R{\mathbb{R}}

\def\D{\mathscr{D}}

\def\L{\mathscr{L}}

\def\phi{\varphi}

\newcommand{\clus}[1][{}]{\displaystyle\mathop{\ifthenelse{\equal{#1}{}}
    {\hbox{\rm clus}}{\hbox{\rm $#1$-clus}}}\limits}
\newcommand{\cl}[2][{}]{\ifthenelse{\equal{#1}{}}
    {\overline{#2}}{\overline{#2}^{#1}}}
\newcommand{\clos}[1][{}]{\mathop{\ifthenelse{\equal{#1}{}}
    {\hbox{\rm cl}}{\hbox{\rm $#1$-cl}}}\nolimits}
\newcommand{\climsup}[1][{}]{\mathop{\ifthenelse{\equal{#1}{}}
    {\hbox{$C$-$\overline{\hbox{\rm Lim}}$}}
    {\hbox{$C$-$\overline{\hbox{\rm Lim}}^{#1}$}}}\limits}
\def\over#1{{\mbox{\tiny $\begin{aligned}#1\end{aligned}$}}}

\def\co{\mathop{\rm co}\nolimits}

\def\clco{\mathop{\overline{\mbox{\rm co}}}\nolimits}

\def\dlim_#1{\lim_{\mbox{\tiny$\begin{array}{c}#1\end{array}$}}}
\def\dlimsup_#1{\limsup_{\mbox{\tiny$\begin{array}{c}#1\end{array}$}}}

\newtheorem{theorem}{Theorem}[section]
\newtheorem*{theorem*}{Theorem}
\def\Thm#1#2{\ifthenelse{\equal{#1}{*}}{\begin{theorem*}#2\end{theorem*}}
             {\begin{theorem}\label{T#1}#2\end{theorem}}}
\newtheorem{Atheorem}{Theorem}

\newtheorem{proposition}[theorem]{Proposition}
\newtheorem*{proposition*}{Proposition}
\def\Prp#1#2{\ifthenelse{\equal{#1}{*}}{\begin{proposition*}#2\end{proposition*}}
             {\begin{proposition}\label{P#1}#2\end{proposition}}}

\newtheorem{corollary}[theorem]{Corollary}
\newtheorem*{corollary*}{Corollary}
\def\Cor#1#2{\ifthenelse{\equal{#1}{*}}{\begin{corollary*}#2\end{corollary*}}
             {\begin{corollary}\label{C#1}#2\end{corollary}}}

\newtheorem{lemma}[theorem]{Lemma}
\newtheorem*{lemma*}{Lemma}
\def\Lem#1#2{\ifthenelse{\equal{#1}{*}}{\begin{lemma*}#2\end{lemma*}}
             {\begin{lemma}\label{L#1}#2\end{lemma}}}

\newtheorem{remark}[theorem]{Remark}
\newtheorem*{remark*}{Remark}
\def\Rem#1#2{\ifthenelse{\equal{#1}{*}}{\begin{remark*}#2\end{remark*}}
             {\begin{remark}\label{R#1}#2\end{remark}}}

\newtheorem{example}[theorem]{Example}

\newtheorem*{example*}{Example}
\def\Exa#1#2{\ifthenelse{\equal{#1}{*}}{\begin{example*}\rm #2\end{example*}}
             {\begin{example}\label{Ex#1}\rm #2\end{example}}}

\def\eq#1{{\rm(\ref{E#1})}}

\def\Eq#1#2{\ifthenelse{\equal{#1}{*}}
  {\begin{equation*}\begin{aligned}#2\end{aligned}\end{equation*}}
  {\begin{equation}\begin{aligned}\label{E#1}#2\end{aligned}\end{equation}}}

\begin{document}
\begin{flushright}
\textit{Submitted to:} Annals of Mathematics 
\end{flushright}
\vspace{5mm}

\date{\today}

\title[Infinite Dimensional Generalized Jacobian]
{Generalized Jacobian for Functions with Infinite
Dimensional Range and Domain}

\author[Zs. P\'ales]{Zsolt P\'ales}
\address{Institute of Mathematics, University of Debrecen, 
H-4010 Debrecen, Pf.\ 12, Hungary}
\email{pales@math.klte.hu}

\author[V. Zeidan]{Vera Zeidan}
\address{Department of Mathematics, Michigan State University,
East Lansing, MI 48824, U.S.A.}
\email{zeidan@math.msu.edu}


\thanks{Research of the first author is supported by the Hungarian
Scientific Research Fund (OKTA) under grant K62316.
Research of the second author is supported by the National Science 
Foundation under grant DMS-0306260.}

\begin{abstract}
In this paper, locally Lipschitz functions acting between infinite
dimensional normed spaces are considered. When the range is a dual
space and satisfies the Radon--Nikod\'ym property, Clarke's
generalized Jacobian will be extended to this setting. 
Characterization and fundamental properties of the extended 
generalized Jacobian are established including the nonemptiness, the
$\beta$-compactness, the $\beta$-upper semicontinuity, and a
mean-value theorem. A connection with known notions is provided and
chain rules are proved using key results developed. This included the
vectorization and restriction theorem, and the extension theorem. 
Therefore, the generalized Jacobian introduced in this paper is proved
to enjoy all the properties required of a derivative like-set.
\end{abstract}

\maketitle

\section{\bf Introduction}

The subject of nonsmooth analysis focuses on the study of a
derivative-like object for nonsmooth functions. When the function is a
convex real-valued, the notion of subgradient was introduced in the 
late fifties by Rockafellar in \cite{Roc70}, and in 
the references therein. Since then, the focus shifted to finding
derivative-like objects for nonconvex, in particular, locally Lipschitz 
functions acting between two normed spaces $X$ and $Y$. 
%It is worth noting that there is another direction 
%f research pertaining the class of locally Lipschitz functions, namely, 
%the study of their differentiability properties (cf.\ \cite{Aro76}, 
%\cite{Chr73}, \cite{CPT05}, \cite{JLPS02a}, \cite{LP00}, \cite{Man73}, 
%and \cite{Phe78}). 

When $X$ and $Y$ are {\it both} finite dimensional normed spaces and
$f:\D\to Y$ is a {\it vector-valued} locally Lipschitz function,
Clarke introduced in \cite{Cla76c}, \cite{Cla83} the notion of the
{\it generalized Jacobian} based on Rademacher's theorem which gives
the almost everywhere differentiability of locally Lipschitz
functions. This generalized Jacobian is  
\Eq{CJ}{
  \partial f(p)
   :=\co\big\{A\in\L(X,Y)\mid\exists\,(x_i)_{i\in\N}\mbox{ in }\Omega(f)
        : \lim_{i\to\infty} x_i= p\, \mbox{ and } 
          \lim_{i\to\infty} Df(x_i)= A\big\}, 
}
where $\Omega(f)$ denotes the set of the points of $\D$ where $f$ is
differentiable, it is of full Lebesgue measure.

Another but related generalized Jacobian based also on the use of
Rademacher's theorem was proposed for the same setting by Pourciau in
\cite{Pou77} and was defined as 
\Eq{PJ}{
  \partial^P f(p):=\bigcap_{\delta>0} \clco
     \big\{D f(x):x\in (p+\delta B_X)\cap\Omega(f)\big\}.
}
One can show that Clarke's and Pourciau's generalized Jacobians are
equivalent. These two objects are nonempty, 
due to Rademacher's
theorem. Furthermore, in terms of these Jacobians, results have been
derived pertaining optimality conditions, implicit functions theorems,
metric regularity, and calculus rules including the sum rule and the
chain rule. Thereby, it has already been shown that these generalized 
Jacobians are successful approximations of $f$ by linear operators. 
 
When establishing calculus rules such as the sum and/or the chain
rules, a fundamental property, namely, the ``blindness'' of the
generalized Jacobian with respect to sets of Lebesgue measure zero, or
null sets is needed. The blindness of Clarke generalized gradient was
established by Clarke in \cite{Cla83} and that of Clarke's generalized
Jacobian was shown by Warga \cite{War81a} and by Fabian \& Preiss
\cite{FP87a}. 

Thibault in \cite{Thi82a} extended Clarke's notion of generalized
Jacobian, equation \eq{CJ}, to the case where $X$ and $Y$ are
infinite dimensional separable Banach spaces such that $Y$ is
{\it reflexive}. This extension was based on the dense differentiability
of locally Lipschitz functions (cf.\ Aronszajn \cite{Aro76},
Christensen \cite{Chr73}, and Phelps \cite{Phe78}). Thibault's
definition is 
\Eq{TJ}{
  \partial_{H} f(p)
   :=\clco\big\{A\in\L(X,Y)\mid\exists\,(x_i)_{i\in\N}\mbox{ in }H
        : \lim_{i\to\infty} x_i= p\, \mbox{ and } 
          \lim_{i\to\infty} Df(x_i)= A\big\}, 
}
where $H$ is a subset of $\D$ on which $f$ is G\^ateaux-differentiable
and such that $\D\setminus H$ is a Haar-null set in $\D$. The notion in
\eq{TJ} depends on the choice of the set $H$, and hence, unlike Clarke's 
generalized Jacobian, is not known to be blind with respect to the
Haar-null sets. In other words, the notion in \eq{TJ}, assigns to
every locally Lipschitz function, not a single object but rather a
family of generalized Jacobians that is parametrized by certain null
sets. Thus, based on this approach all the chain rules derived in 
\cite{Thi82a} are in terms of the Haar null set $H$.

Other notions are known in the infinite dimensional setting, such as 
the notion of derivate containers in \cite{War76b}, \cite{War81a}; 
the concepts of screens and ``fans'' in \cite{Hal76a}, \cite{Hal76b}; 
the concept of shields \cite{Swe77}; Ioffe's fan derivative
\cite{Iof81c}, and the notion of coderivatives developed in
\cite{MS96c}. Most of these notions are not given in terms of relevant
sets of {\it linear operators}. A relatively recent survey on the
different subdifferentials and their properties is given in
\cite{BZ99a} where also an extended list of references could be found. 

In recent papers \cite{PZ*a}and \cite{PZ*b} Clarke's generalized
Jacobian \eq{CJ} was extended to the case when $X$ was any {\it normed}
space and $Y$ was a finite dimensional space. In these references the
generalized Jacobian was defined to be a {\it set of linear operators}
from $X$ to $Y$. When the domain is infinite dimensional and the image
space is $\R$, the notion introduced in \cite{PZ*a} and \cite{PZ*b}
coincides with Clarke's generalized gradient which is defined as 

\Eq{Cg}{
  \partial^c f(p)
   := \big\{\zeta\in X^*\mid\langle\zeta, h\rangle\leq f^{\circ}(p, h), 
       \,\, \forall h\in X\big\},
}   
where   
\Eq{Cdd}{
  f^{\circ}(p,h):=\limsup_{\over{x&\to p\\[-2mm]t&\to 0^+}}
      \frac{f(x+th)-f(x)}{t}
}
is Clarke's {\it generalized directional derivative}.

In \cite{PZ*a}and \cite{PZ*b}, the nonemptiness, the $w^*$-compactness, the
convexity, and the upper semicontinuity property of the extended generalized 
Jacobian were derived. Furthermore,  a chain rule for
the composition of nonsmooth locally Lipschitz maps with finite
dimensional ranges was established. 

The difficulty caused by the infinite dimensionality of the domain was
handled in \cite{PZ*a}and \cite{PZ*b} by introducing an intermediate
Jacobian $\partial_L f$ defined on finite dimensional spaces $L$ so
that Rademacher theorem remains applicable. 

In this paper we are interested in extending the definition of Clarke's
generalized 
Jacobian to the case when in addition to the domain also the range is
infinite dimensional. In this case, two extra difficulties manifest. 
The first is the differentiability issue related to the Rademacher
theorem in infinite dimension. This issue will be handled by taking 
image spaces satisfying the Radon--Nikod\'ym property. This implies that 
the restriction of a Lipschitz function  $f:\D\to Y$ to a 
{\it finite dimensional} domain is almost everywhere differentiable 
(cf.\ \cite{BL00}[Prop. 6.41]).

The second difficulty is pertaining finding a topology in the space of 
linear operators $\L(X,Y)$, where the generalized Jacobian lives, so that
bounded sequences would have cluster points in this topology. To
overcome this difficulty, we also assume that the image space $Y$ is a
dual of a normed space. 

The goal of this paper is to provide a generalized Jacobian for
locally Lipschitz functions defined between infinite dimensional
normed spaces with the range $Y$ is a dual space and satisfies the
Radon--Nikod\'ym property. We shall show that our generalized Jacobian
enjoys all the fundamental properties desired from a derivative set.

In Section 2 we introduce the $\beta$-topology on the space of linear
operators $\L(X,Y)$. This is a $w^*$-operator topology induced by the
predual of $Y$. We prove an analog of the Banach--Alaoglu theorem
 as well as an extension theorem  which will be
crucial for the proof of the nonemptiness of the generalized Jacobian. 
In this section we also derive results related to various upper
semicontinuity properties which will be repeatedly used in the
subsequent section. In Section 3 the $L$-Jacobian, $\partial _Lf(p)$,
and the generalized Jacobian are defined as an extension of Pourciau's
notion, equation \eq{PJ}, to infinite dimensional spaces. We also show
that our generalized Jacobian could be equivalently defined in terms
of cluster points, which is a definition  that corresponds to Clarke's
approach. Basic properties and a characterization of the generalized 
Jacobian are established. A main tool named the ``restriction and the
vectorization'' theorem is developed which is central for
deriving many results in the rest of the paper.

Relationships to Thibault's limit set, to a Ioffe type fan derivative and 
to Mordukhovich coderivative are given in Section 4. A generalization of
Lebourg mean-value theorem is obtained as well as that the
generalized Jacobian is a $w^*$-Hadamard prederivative. We also
characterize the cases when the generalized Jacobian is a strict
norm-G\^ateaux or a strict $w^*$-Fr\'echet prederivative. In Section 5 we
derive two chain rules: a nonsmooth-smooth, and a nonsmooth-nonsmooth
one. Their proofs evoke most of the results and properties
established in the previous sections. As a consequence, a sum rule
follows. Finally, in Section 6 we develop results for the generalized
Jacobian of continuous selections.

\def\accepted{accepted for publications}

\def\DL#1{#1}\def\sort#1{}\def\cprime{$'$} \def\Dbar{\leavevmode\lower.6ex\hbox
  to 0pt{\hskip-.23ex \accent"16\hss}D}
  \def\cfac#1{\ifmmode\setbox7\hbox{$\accent"5E#1$}\else
  \setbox7\hbox{\accent"5E#1}\penalty 10000\relax\fi\raise 1\ht7
  \hbox{\lower1.15ex\hbox to 1\wd7{\hss\accent"13\hss}}\penalty 10000
  \hskip-1\wd7\penalty 10000\box7}
  \def\cftil#1{\ifmmode\setbox7\hbox{$\accent"5E#1$}\else
  \setbox7\hbox{\accent"5E#1}\penalty 10000\relax\fi\raise 1\ht7
  \hbox{\lower1.15ex\hbox to 1\wd7{\hss\accent"7E\hss}}\penalty 10000
  \hskip-1\wd7\penalty 10000\box7}
  \def\cfgrv#1{\ifmmode\setbox7\hbox{$\accent"5E#1$}\else
  \setbox7\hbox{\accent"5E#1}\penalty 10000\relax\fi\raise 1\ht7
  \hbox{\lower1.05ex\hbox to 1\wd7{\hss\accent"12\hss}}\penalty 10000
  \hskip-1\wd7\penalty 10000\box7}
  \def\cfudot#1{\ifmmode\setbox7\hbox{$\accent"5E#1$}\else
  \setbox7\hbox{\accent"5E#1}\penalty 10000\relax\fi\raise 1\ht7
  \hbox{\raise.1ex\hbox to 1\wd7{\hss.\hss}}\penalty 10000 \hskip-1\wd7\penalty
  10000\box7} \def\cprime{$'$}
  \def\ocirc#1{\ifmmode\setbox0=\hbox{$#1$}\dimen0=\ht0 \advance\dimen0
  by1pt\rlap{\hbox to\wd0{\hss\raise\dimen0
  \hbox{\hskip.2em$\scriptscriptstyle\circ$}\hss}}#1\else {\accent"17 #1}\fi}
\providecommand{\bysame}{\leavevmode\hbox to3em{\hrulefill}\thinspace}
\providecommand{\MR}{\relax\ifhmode\unskip\space\fi MR }
\providecommand{\MRhref}[2]{%
  \href{http://www.ams.org/mathscinet-getitem?mr=#1}{#2}
}
\providecommand{\href}[2]{#2}

\end{document}